\documentclass[fleqn]{article}

\newcommand{\R}{\mathbb{R}}

\newcommand{\bq}{\mbox{{\boldmath $\mathcal{Q}$}}}
\newcommand{\bp}{\mbox{{\boldmath $\mathcal{P}$}}}

\newtheorem{thm}{Theorem}[section]
\newtheorem{lem}[thm]{Lemma}
\newtheorem{prop}[thm]{Proposition}

\usepackage{amsmath}
\usepackage{amssymb}

\begin{document}
\title{Approximate Nonnegative Matrix Factorization via Alternating Minimization}
\author{Lorenzo Finesso \\ ISIB--CNR \\ Corso Stati Uniti, 4 \\ 35127 Padova -- Italy \\ {\tt finesso@isib.cnr.it}
\and Peter Spreij \\ Korteweg-de Vries Institute for
Mathematics \\
Universiteit van Amsterdam \\
Plantage Muidergracht 24 \\
1018 TV Amsterdam -- The Netherlands \\
{\tt spreij@science.uva.nl}
}
\maketitle
\begin{abstract}
\noindent In this paper we consider the Nonnegative Matrix
Factorization (NMF) problem: given an (elementwise) nonnegative
matrix $V \in \R_+^{m\times n}$ find, for assigned $k$,
nonnegative matrices $W\in\R_+^{m\times k}$ and $H\in\R_+^{k\times
n}$ such that $V=WH$. Exact, non trivial, nonnegative
factorizations do not always exist, hence it is interesting to
pose the approximate NMF problem. The criterion which is commonly
employed is I-divergence between nonnegative matrices. The problem
becomes that of finding, for assigned $k$, the factorization $WH$
closest to $V$ in I-divergence. An i\-te\-ra\-ti\-ve algorithm, EM
like, for the construction of the best pair $(W, H)$ has been
proposed in the literature. In this paper we interpret the
algorithm as an alternating minimization procedure \`a la
Csisz\'ar-Tusn\'ady and investigate some of its stability
properties. NMF is widespreading as a data analysis method in
applications for which the positivity constraint is relevant.
There are other data analysis methods which impose some form of
nonnegativity: we discuss here the connections between NMF and
Archetypal Analysis. An interesting system theoretic application
of NMF is to the problem of approximate realization of Hidden
Markov Models.

\end{abstract}

\newpage

\section{Introduction}
The approximate Nonnegative Matrix Factorization (NMF) of
nonnegative matrices is a data analysis technique only recently
introduced \cite{leeseung1999, sullivan2000}. Roughly speaking the
problem is to find, for a given nonnegative matrix $V \in
\R_+^{m\times n}$, and an assigned $k$, a pair of nonnegative
matrices $W\in\R_+^{m\times k}$ and $H\in\R_+^{k\times n}$ such
that, in an appropriate sense, $V \approx WH$. EM like algorithms
for the construction of a factorization have been proposed in
\cite{leeseung1999, leeseung2001}. In \cite{sullivan2000} the
connection of these algorithms with the classic alternating
minimization of the I-divergence \cite{ct1984} has been pointed
out but not fully investigated. In this paper we pose the NMF
problem as a minimum I-divergence problem that can be solved by
alternating minimization and derive, from this point of view, the
algorithm proposed in \cite{leeseung1999}.

Although only recently introduced the NMF has found many
applications as a data reduction procedure and has been advocated
as an alternative to Principal Components Analysis (PCA) in cases
where the positivity constraint is relevant (typically image
analysis). The title of \cite{sullivan2000} is a clear indication
of this point of view, but a complete analysis of the relations
between NMF and PCA is still lacking. Other data analysis methods
proposed in the literature enforce some form of positivity
constraint and it is useful to investigate the connection between
NMF and these methods. An interesting example is the so called
Archetypal Analysis (AA) technique \cite{cb1994}. Assigned a
matrix $X \in \R^{m\times n}$ and an integer $k$, the AA problem
is to find, in the convex hull of the columns of $X$, a set of $k$
vectors whose convex combinations can optimally represent $X$. To
understand the relation between NMF and AA we choose the $L_2$
criterion for both problems. For any matrix $A$ and positive
definite matrix $\Sigma$ define $||A||_\Sigma = (\rm{tr} (A^T 
\Sigma
A))^{1/2} $. Denote $||A||_I = ||A||$. The solution of the NMF problem is
then
$$
(W, H) = \arg \min_{W, H} || V - WH ||
$$
where the minimization is constrained to the proper set of
matrices. The solution to the AA problem is given by the pair of
column stochastic matrices $(A, B)$ of respective sizes $k \times
n$ and $m \times k$ such that $||X - XBA||$ is minimized (the
constraint to column stochastic matrices is imposed by the
convexity). Since $||X - XBA|| = ||I - BA||_{X^TX}$ the solution
of the AA problem is
$$
(A, B) = \arg \min_{A, B} || I - BA ||_{X^TX}.
$$
AA and NMF can therefore be viewed as special cases of a more
general problem which can be stated as follows. Given any matrix
$P \in \R_+^{m\times n}$, any positive definite matrix $\Sigma$,
and any integer $k$, find the best nonnegative factorization $P
\approx Q_1Q_2$ (with $Q_1 \in\R_+^{m\times k}, \,\, Q_2
\in\R_+^{k\times n}$) in the $L_2$ sense, {\it i.e.}
$$
(Q_1, Q_2) =  \arg \min_{Q_1, Q_2} ||P - Q_1Q_2||_\Sigma.
$$
Our interest in NMF stems from the system theoretic problem of
approximate realization (or order reduction) of Hidden Markov
Models. Partial results have already been obtained \cite{fs2002}.

\section{Preliminaries and problem statement}\label{sec:problem}

The NMF is a long standing problem in linear algebra \cite{haz,
phs}. It can be stated as follows. Given $V \in \R_+^{m\times n}$,
and $1 \le k \le \min(m,n)$, find a pair of matrices
$W\in\R_+^{m\times k}$ and $H\in\R_+^{k\times n}$ such that $V =
WH$.  The smallest $k$ for which a factorization exists is called
the positive rank of $V$, denoted ${\rm prank}(V)$. This
definition implies that ${\rm rank}(V) \le {\rm prank}(V) \le
\min(m,n)$. It is well known that ${\rm prank}(V)$ can assume all
intermediate values, depending on $V$. Examples for which
nonnegative factorizations do not exist, and examples for which
factorization is possible only for $k > {\rm rank}(V)$ are easily
constructed \cite{haz}. The ${\rm prank}$ has been characterized
only for special classes of matrices \cite{phs} and algorithms for
the construction of a NMF are not known. The approximate NMF has
been recently introduced in \cite{leeseung1999} independently from
the exact NMF problem. The set-up is the same, but instead of
exact factorization it is required that $V \approx WH$ in an
appropriate sense.
In \cite{leeseung1999} and in this paper the approximation is to
be understood in the sense of minimum I-divergence. For two
nonnegative matrices (or vectors) $M=(M_{ij})$ and $N=(N_{ij})$ of
the same size the I-divergence is defined as
\begin{equation*}
D(M||N)= \sum_{ij}
(M_{ij}\log\frac{M_{ij}}{N_{ij}}-M_{ij}+N_{ij}),
\end{equation*}
with the conventions $0/0=0$, $0\log 0=0$ and $p/0=\infty$ for
$p>0$.
>From the inequality $x\log x\geq x-1$ it follows that $D(M||N)\geq
0$ with equality iff $M=N$. The problem of approximate NMF is to
find
\begin{equation*} \arg \min_{W,H}D(V||WH).
\end{equation*}
It can be shown that, if $V_{ij}>0$, the minimum is attained.
Dropping constants the problem is equivalent to finding
\begin{equation*}
\max_{W,H} F(W,H):= \sum_{ij} (V_{ij}\log (WH)_{ij}-(WH)_{ij}).
\end{equation*}
Clearly the solution is not unique. In order to rule out too many
trivial multiple solutions, we impose the condition that $H$ is
row stochastic, so $\sum_jH_{lj}=1$ for all $l$. This is not a
restriction. Indeed, excluding without loss of generality the case
where $H$ has one or more zero rows, let $h$ be the diagonal
matrix with elements $h_i=\sum_j H_{ij}$, then
$WH=\tilde{W}\tilde{H}$ with $\tilde{W}=Wh$, $\tilde{H}=h^{-1}H$
and $\tilde{H}$ is by construction row stochastic. The convention
that $H$ is row stochastic still doesn't rule out non-uniqueness.
Think e.g.\ of post-multiplying $W$ with a permutation matrix
$\Pi$ and pre-multiplying $H$ with $\Pi^{-1}$.

Although the function $F$ is concave in each of its arguments $W$
and $H$ separately, it does not have this property as a function
of two variables. Hence $F$ may have several (local) maxima, that
may prevent numerical algorithms for a global maximum search to
converge to the global maximizer.

Let $e$ ($e^\top$) be a column (row) vector of appropriate
dimension whose elements are all equal to one. The (constrained)
problem we will look at is then
\begin{equation}\label{maxc}
\max_{W,H: He=e} F(W,H).
\end{equation}
Notice that the constrained problem~(\ref{maxc}) can be rewritten
as
\begin{equation*}
\max_{W,H: He=e} F(W,H):= \sum_{ij} (V_{ij}\log (WH)_{ij}-W_{ij}).
\end{equation*}
To carry out the maximization numerically
~\cite{leeseung1999,leeseung2001} propose an iterative algorithm.
Denoting by $W^n$ and $H^n$ the matrices at step $n$, the update
equations are the following
\begin{align}
W^{n+1}_{il} & =  \sum_j V_{ij}\frac{W^n_{il}H^n_{lj}}{(W^nH^n)_{ij}}\label{eq:wbar2}\\
H^{n+1}_{lj} & =  \sum_i V_{ij}\frac{W^n_{il}H^n_{lj}}{(W^nH^n)_{ij}}/ \sum_i
\sum_j V_{ij}\frac{W^n_{il}H^n_{lj}}{(W^nH^n)_{ij}}.\label{eq:hbar2}
\end{align}
There is no rationale for this algorithm although the update
steps~(\ref{eq:wbar2}) and~(\ref{eq:hbar2}) are like those in the
EM algorithm, known from statistics, see~\cite{em}. Likewise the
convergence properties of the algorithm are unclear. In the next
section we will cast the maximization problem in a different way
that provides more insight in the specific form of the update
equations.

\section{Lifted version of the problem}

In this section we lift the I-divergence minimization problem to
an equivalent minimization problem where the `matrices' (we should
speak of {\em tensors}) have three indices. Because we insist on
probabilistic interpretations we change notations as follows.
$P \in \R_+^{m\times n}$ is a given, fixed matrix and 
\bigskip\\ $\bp=\{\mathbf{P}\in\mathbb{R}^{m\times k\times n}_+:
\sum_l\mathbf{P}(ilj)=P(ij)\}$, 
\bigskip\\ $ \bq = \{\mathbf{Q}\in\mathbb{R}^{m\times k\times n}_+:
\mathbf{Q}(ilj)=Q_-(il)Q_+(lj), \,\,\, Q_-(il),\,\, Q_+(lj)\ge0,
\,\, Q_+e=e\}$,
\bigskip\\ $ \mathcal{Q} = \{ Q\in \mathbb{R}^{m\times n}_+:
Q(ij)=\sum_l\mathbf{Q}(ilj) \quad {\rm for \,\,  some} \,\,
\mathbf{Q}\in \bq \}. $
\bigskip\\
Notice that $\mathcal{Q}$ is the class of $m\times n$ matrices
that admit exact NMF of size $k$. In the notation of section $2$,
$V$ has become $P$, and $W, H$ are now $Q_-, Q_+$ respectively.
\medskip\\
The following observation (whose proof is elementary,
see~\cite{ps}) motivates our approach.

\begin{lem}\label{lemma:61}
$P$ can be factorized as $P=Q_-Q_+$ iff $\bp\cap\bq\neq\emptyset$,
so iff there exists a $\mathbf{P}\in\bp$ and $\mathbf{Q}\in\bq$
such that $\mathbf{P}=\mathbf{Q}$.
\end{lem}
For a probabilistic interpretation of this lemma, and of the
results below, we assume (without loss of generality) that
$\mathbf{P}$ represents the joint distribution of a three
dimensional random vector. Suppose that $Y_-$ and $Y_+$ are finite
valued random variables defined on a probability space
$(\Omega,\mathcal{F},\mathbb{P})$ whose joint distribution is
given by $\mathbb{P}(Y_-=i,Y_+=j)=P(ij)$. Then the content of the
lemma is that there exists a finite valued random variable $X$
such that $Y_-$ and $Y_+$ are conditionally independent given $X$
iff $P=Q_-Q_+$. The matrix $Q_-$ then gives the joint distribution
of $Y_-$ and $X$ by $Q_-(il)=\mathbb{P}(Y_-=i,X=l)$, whereas the
matrix $Q_+$ can be interpreted as conditional distributions of
$Y^+$ given $X$ via $Q_+(lj)=\mathbb{P}(Y_+=j|X=l)$. Moreover, in
this case we have $\mathbb{P}(Y_-=i,X=l,Y_+=j)=\mathbf{Q}(ilj)$.
To see this we write the conditional independence relation
\[
\mathbb{P}(Y_-=i,Y_+=j|X=l)=\mathbb{P}(Y_-=i|X=l)\mathbb{P}(Y_+=j|X=l)
\]
in equivalent form as
\[
\mathbb{P}(Y_-=i,X=l,Y_+=j)=\mathbb{P}(Y_-=i,X=l)\mathbb{P}(Y_+=j|X=l),
\]
from which the above statements immediately follow.

\section{Two partial minimization problems}

In this section we consider the following two minimization
problems. In the first one we minimize for given
$\mathbf{Q}\in\bq$ the I-divergence $D(\mathbf{P}||\mathbf{Q})$
over $\mathbf{P}\in\bp$. In the second problem we minimize for
given $\mathbf{P} \in\bp$ the I-divergence
$D(\mathbf{P}||\mathbf{Q})$ over $\mathbf{Q}\in\bq$. The unique
solution $\mathbf{P}^*=\mathbf{P}^*(\mathbf{Q})$ to the first
problem can be computed analytically and is given by
\begin{equation}\label{eq:p*}
\mathbf{P}^*(ilj)=\frac{\mathbf{Q}(ilj)P(ij)}{Q(ij)},
\end{equation}
where $Q(ij)=\sum_l\mathbf{Q}(ilj)$. A direct computation gives
the useful relation
\begin{equation*}
D(\mathbf{P}^*(\mathbf{Q})||\mathbf{Q})=D(P||Q).
\end{equation*}
\medskip\\
The interpretation in terms of random variables is that for a
given probability measure $\mathbb{Q}$, random variables
$Y_-,X,Y_+$ with law
$\mathbb{Q}(Y_-=i,X=l,Y_+=j)=\mathbf{Q}(ilj)$, the best
approximating model $\mathbb{P}^*$ with marginal distribution of
$Y=(Y_-,Y_+)$  described by $P$ is given by
\begin{align*}
\mathbf{P}^*(ilj) & = \mathbb{P}(Y_-=i,X=l,Y_+=j) \\
& = \mathbb{Q}(X=l|Y_-=i,Y_+=j)P(i,j).
\end{align*}
Equivalently, we can say that $\mathbb{P}^*$ is such that the marginal
distribution of $Y$ under $\mathbb{P}^*$ is given by $P$ and the conditional
distribution of $X$ given $Y$ under $\mathbb{P}^*$ is equal to the
conditional distribution under $\mathbb{Q}$. Below we will see that
this is not a coincidence.
\medskip\\
The solution $\mathbf{Q}^*=\mathbf{Q}^*(\mathbf{P})$ to the second
problem is given by
\begin{align}
Q^*_-(il) & = \sum_j\mathbf{P}(ilj)\label{eq:q-} \\
Q^*_+(lj) & = \frac{\sum_i \mathbf{P}(ilj)}{\sum_{ij}\mathbf{P}(ilj)}.\label{eq:q+}
\end{align}
The interpretation in probabilistic terms is that for a given
distribution $\mathbb{P}$ of $(Y_-,X,Y_+)$, the best model $\mathbb{Q}^*$
that makes $Y_-$ and $Y_+$ conditionally independent given $X$ is such
that
\[
\mathbb{Q}^*(Y_-=i,X=l)=\mathbb{P}(Y_-=i,X=l)
\]
and
\[
\mathbb{Q}^*(Y_+=j|X=l,Y_-=i)=\mathbb{Q}^*(Y_+=j|X=l)=\mathbb{P}(Y_+=j|X=l).
\]
We see that the optimal solution $\mathbb{Q}^*$ is such that the
marginal distributions of $(X,Y_-)$ under $\mathbb{P}$ and
$\mathbb{Q}^*$ coincide and that the same happens for the conditional
distributions of $Y_+$ given $X$. Again, this is not a coincidence, as
we will explain below. First we will state
for the two partial minimization problems above  the following two Pythagorean rules.
\begin{lem}\label{lemma:pyth}
For fixed $\mathbf{P}$ and $\mathbf{Q}^*=\mathbf{Q}^*(\mathbf{P})$
it holds that for any  $\mathbf{Q} \in\bq$
\begin{equation}\label{eq:pythq}
D(\mathbf{P}||\mathbf{Q})=D(\mathbf{P}||\mathbf{Q}^*)+D(\mathbf{Q}^*||\mathbf{Q}),
\end{equation}
whereas for fixed $\mathbf{Q}$ and
$\mathbf{P}^*=\mathbf{P}^*(\mathbf{Q})$ it holds that for any
$\mathbf{P} \in\bp$
\begin{equation}\label{eq:pythp}
D(\mathbf{P}||\mathbf{Q})=D(\mathbf{P}||\mathbf{P}^*)+D(\mathbf{P}^*||\mathbf{Q}),
\end{equation}
and
\begin{equation}\label{eq:p0q0}
D(\mathbf{P}^*||\mathbf{Q})=D(P||Q),
\end{equation}
where $Q$ is given by $Q(ij)=\sum_l\mathbf{Q}(ilj)$.
\end{lem}
{\bf Proof.}
To prove the first relation we first introduce some notation. Let
$\mathbf{P}(il\cdot)=\sum_j \mathbf{P}(ilj)$, $\mathbf{P}(\cdot lj)=\sum_i\mathbf{P}(ilj)$ and
$\mathbf{P}(j|l)=\mathbf{P}(\cdot lj)/\sum_j\mathbf{P}(\cdot lj)$. For $\mathbf{Q}$ we use similar notation and so we
have
$\mathbf{Q}(il\cdot)=Q_-(il)$ and $\mathbf{Q}(j|l)=Q_+(lj)/\sum_jQ_+(lj)$ and
$Q^*_-(il)=\mathbf{P}(il\cdot)$ and $Q^*_+(lj)=\mathbf{P}(j|l)$.
Consider
\begin{align*}
D(\mathbf{P}||\mathbf{Q})-D(\mathbf{P}||\mathbf{Q}^*) &=\sum_{ilj}\mathbf{P}(ilj)\log
\frac{\mathbf{P}(il\cdot)}{Q_-(ij)}+\log\frac{\mathbf{P}(j|l)}{Q_+(lj)}\\
& = \sum_{il}\mathbf{P}(il\cdot)\log
\frac{\mathbf{P}(il\cdot)}{Q_-(ij)} + \sum_{lj}\mathbf{P}(\cdot lj)\log\frac{\mathbf{P}(j|l)}{Q_+(lj)}.
\end{align*}
On the other hand we have
\begin{align*}
D(\mathbf{Q}^*||\mathbf{Q}) & = \sum_{ilj}
\mathbf{P}(il\cdot)\mathbf{P}(j|l)(\log\frac{\mathbf{P}(il\cdot)}{Q_-(il)}+\log\frac{\mathbf{P}(j|l)}{Q_+(lj)} \\
& = \sum_{il}\mathbf{P}(il\cdot)\log\frac{\mathbf{P}(il\cdot)}{Q_-(il)}+
\sum_{lj}\mathbf{P}(\cdot lj)\log\frac{\mathbf{P}(j|l)}{Q_+(lj)}.
\end{align*}
The first assertion follows.
The second Pythagorean rule follows from
\begin{eqnarray*}
\lefteqn{D(\mathbf{P}||\mathbf{P}^*)+ D(\mathbf{P}^*||\mathbf{Q})}\\
& = &  \sum_{ilj}\mathbf{P}(ilj)\log\frac{\mathbf{P}(ilj)Q(ij)}{\mathbf{Q}(ilj)\mathbf{P}(ij)} +
\sum_{ilj}\mathbf{Q}(ilj)\frac{P(ij)}{Q(ij)}\log \frac{P(ij)}{Q(ij)} \\
& = & \sum_{ilj}\mathbf{P}(ilj)\log\frac{\mathbf{P}(ilj)}{\mathbf{Q}(ilj)} +
\sum_{ilj}\mathbf{P}(ilj)\log\frac{Q(ij)}{P(ij)} \\
& & \mbox{} + \sum_{ij}Q(ij)\frac{P(ij)}{Q(ij)}\log \frac{P(ij)}{Q(ij)} \\
& = &   D(\mathbf{P}||\mathbf{Q}).
\end{eqnarray*}
\hfill$\square$\medskip\\
For a probabilistic interpretation of the $\mathbf{P}^*$ and
$\mathbf{Q}^*$ above as well as the Py\-tha\-go\-rean rules we use
a general result on the I-divergence between two joint laws of a
random vector $(U,V)$. We denote the law of this vector under
probability measures $\mathbb{P}$ and $\mathbb{Q}$ by $P^{U,V}$
and $Q^{U,V}$. The conditional distributions of $U$ given $V$ are
summarized by the matrices $P^{U|V}$ and $Q^{U|V}$, with the
obvious convention $P^{U|V}(ij)=\mathbb{P}(U=i|V=j)$ and likewise
for $Q^{U|V}$.
\begin{lem}
It
holds that
\begin{equation}\label{eq:duv}
D(P^{U,V}||Q^{U,V})=\mathbb{E}_\mathbb{P} D(P^{U|V}||Q^{U|V}) +
D(P^V||Q^V),
\end{equation}
where
\[
D(P^{U|V}||Q^{U|V}) = \sum_i P(U=i|V)\log\frac{P(U=i|V)}{Q(U=i|V)}.
\]
\end{lem}
{\bf Proof.} This follows from elementary manipulations.
\hfill $\square$\medskip\\
The above relation can be refined as follows. Suppose that $V$ is
bivariate,
$V=(V_1,V_2)$ say and that $U$ and $V_2$ are conditionally independent
given $V_1$ under $\mathbb{Q}$, so the conditional distribution of $U$
given $V$ is the same as the conditional distribution of $U$ given
$V_1$ under $\mathbb{Q}$. Then the first term on the right hand side of equation~(\ref{eq:duv}) can be
decomposed
as
\begin{equation}\label{eq:duv1}
\mathbb{E}_\mathbb{P} D(P^{U|V}||Q^{U|V}) =\mathbb{E}_\mathbb{P} D(P^{U|V}||P^{U|V_1}) +
\mathbb{E}_\mathbb{P} D(P^{U|V_1}||Q^{U|V_1}).
\end{equation}
We apply this lemma to the first partial minimization problem
above by an appropriate choice of $U$ and $V$. Since
$D(\mathbf{P}^*||\mathbf{Q})=D(P||Q)=D(P^Y||Q^Y)$, where $P^Y$ is
given by $P$, we see that for $U=X$, $V=Y=(Y_-,Y_+)$ the
decomposition~(\ref{eq:pythp}) can alternatively be written as
$\mathbb{E}_\mathbb{P}D(P^{X|Y}||Q^{X|Y}) + D(P||Q)$. Minimizing
$D(\mathbf{P}||\mathbf{Q})$ w.r.t.\ $\mathbf{P}$ under the
condition that the marginal of $\mathbf{P}$ is given by $P$ is
thus equivalent to minimizing the I-divergence between the
conditional distributions $P^{X|Y}$ and $Q^{X|Y}$, and this
clearly happens for
$P^{X|Y}=Q^{X|Y}$.\\
The interpretation of~(\ref{eq:pythq}) is less straightforward.
However, refining~(\ref{eq:pythq}), we have parallel
to~(\ref{eq:duv1})
\begin{align*}
D(\mathbf{P}||\mathbf{Q})= & \mathbb{E}_\mathbb{P}D(P^{Y_+|X,Y_-}||P^{Y_+|X})\\
& \mbox{} +D(P^{Y_-,X}||Q^{Y_-,X})+\mathbb{E}_\mathbb{P}D(P^{Y_+|X}||Q^{Y_+|X}).
\end{align*}
Hence the minimization problem here is to minimize the
I-divergence between the distributions of $(Y_-,X)$ under
$\mathbb{P}$ and $\mathbb{Q}$ and the I-divergence between the
conditional probability measures $P^{Y_+|X}$ and $Q^{Y+|X}$. This
explains the form of the optimal solution
$\mathbf{Q}^*(\mathbf{P})$.
\medskip\\
The next proposition shows that the original minimization of $D(P||Q)$
over nonnegative matrices $Q$ for a given nonnegative matrix $P$ is
equivalent to a double minimization over the sets $\bp$ and $\bq$.
\begin{prop}\label{prop:pqq}
Let $P$ be given.
It holds that
\begin{equation*}
\min_{Q\in\mathcal{Q}}D(P||Q)=\min_{\mathbf{P}\in\bp,\mathbf{Q}\in\bq}D(\mathbf{P}||\mathbf{Q}).
\end{equation*}
\end{prop}
{\bf Proof.} With $\mathbf{P}^*=\mathbf{P}^*(Q)$, the optimal solution of the partial minimization over $\bp$,
we have
\begin{align*}
D(\mathbf{P}||\mathbf{Q})& \geq D(\mathbf{P}^*||\mathbf{Q}) \\
& = D(P||Q) \\
& \geq  \min_{Q\in\mathcal{Q}}D(P||Q).
\end{align*}
It follows that
$\min_{\mathbf{P}\in\bp,\mathbf{Q}\in\bq}D(\mathbf{P}||\mathbf{Q})\geq\min_{Q\in\mathcal{Q}}D(P||Q)$.
\\
Conversely, let $Q^*\in\mathcal{Q}$ be the minimizer of $D(P||Q)$ and let
$\mathbf{Q}$ be a corresponding element in $\bq$. Furthermore, let $\mathbf{P}\in\bp$ be arbitrary. Then we have
\begin{align*}
D(P||Q^*) & \geq D(\mathbf{P}^*(\mathbf{Q})|| \mathbf{Q})\\
& \geq
\min_{\mathbf{P}\in\bp,\mathbf{Q}\in\bq}D(\mathbf{P}||\mathbf{Q}),
\end{align*}
which shows the other inequality.
\hfill$\square$

\section{Alternating minimization algorithm}

The results of the previous section are aimed at setting up an
alternating minimization algorithm for obtaining $\min_Q D(P||Q)$,
where $P$ is a given nonnegative matrix. In view of
proposition~\ref{prop:pqq} we can lift this problem to the $(\bp,
\bq)$ space. Starting with an arbitrary $\mathbf{Q}_1\in\bq$ with
strictly positive elements, we adopt the following recursive
scheme
\begin{equation}\label{eq:altalgo}
\mathbf{Q}_n \to \mathbf{P}_n\to \mathbf{Q}_{n+1}\to \mathbf{P}_{n+1},
\end{equation}
where $\mathbf{P}_n=\mathbf{P}^*(\mathbf{Q}_n)$, $\mathbf{Q}_{n+1}=\mathbf{Q}^*(\mathbf{P}_n)$ and $\mathbf{P}_{n+1}=\mathbf{P}^*(\mathbf{Q}_{n+1})$.
\medskip
The two Pythagorean rules from lemma~\ref{lemma:pyth} now take the
forms
\begin{align*}
D(\mathbf{P}_{n}||\mathbf{Q}_{n+1}) & = D(\mathbf{P}_{n}||\mathbf{P}_{n+1})+D(\mathbf{P}_{n+1}||\mathbf{Q}_{n+1}) \\
D(\mathbf{P}_n||\mathbf{Q}_n) & = D(\mathbf{P}_n||\mathbf{Q}_{n+1})+D(\mathbf{Q}_{n+1}||\mathbf{Q}_n).
\end{align*}
Addition  of these two equations
results in
\begin{align*}
D(\mathbf{P}_{n}||\mathbf{Q}_n) & = D(\mathbf{P}_{n}||\mathbf{P}_{n+1})+D(\mathbf{P}_{n+1}||\mathbf{Q}_{n+1})+D(\mathbf{Q}_{n+1}||\mathbf{Q}_n),
\end{align*}
and together with~(\ref{eq:p0q0}) this becomes
\begin{equation}\label{eq:gain}
D(P||Q_{n})  = D(\mathbf{P}_{n}||\mathbf{P}_{n+1})+D(P||Q_{n+1})+D(\mathbf{Q}_{n+1}||\mathbf{Q}_n).
\end{equation}
This equation also shows that $D(P||Q_{n})  \geq
D(P||Q_{n+1})$.
The procedure outlined in equation~(\ref{eq:altalgo}) will be made
explicit, using equations~(\ref{eq:p*}), (\ref{eq:q+}) and
(\ref{eq:q-}). Since it is our aim to apply the above results to the
problem as sketched in section~\ref{sec:problem}, we now turn back to
the notation of that section. So, instead of $Q_-$ we write $W$,
instead of $Q_+$ we write $H$, instead of $Q$ we write $WH$,
of course these will be endowed with superscript indices $n$ and $n+1$ below,
and $P$ becomes $V$ again.
From~(\ref{eq:altalgo}) we get
$\mathbf{Q}_{n+1}=\mathbf{Q}^*(\mathbf{P}^*(\mathbf{Q}_n))$ and
combining this with the substitution of~(\ref{eq:p*})
into~(\ref{eq:q-}) we obtain--in the original notation--
\[
W^{n+1}_{il}=\sum_j \frac{W^n_{il}H^n_{lj}V_{ij}}{(W^nH^n)_{ij}},
\]
which is just~(\ref{eq:wbar2}). Of course~(\ref{eq:hbar2}) can be
derived similarly.

\section{Discussion of the algorithm}

In the previous section we have shown that the update
rules~(\ref{eq:wbar2}) and~(\ref{eq:hbar2}) are the result of an
alternating minimization procedure. The convergence properties of
the algorithm can be studied using the general results
of~\cite{ct1984}. Due to the similarity with the EM algorithm one
may expect similar convergence properties, see~\cite{wu}.
\medskip\\
At each iteration the I-divergence between $V$ and the $W^nH^n$ is
reduced, equivalently the sequence $F(W^n,H^n)$ is increasing.
This follows from equation~(\ref{eq:gain}). Secondly, once the
algorithm reaches a stationary $(W,H)$-point of $F$ (the partial
derivatives vanish here), the updated values are exactly equal to
the given values. This can be immediately seen by computing the
fist order necessary conditions  for a stationary point and
comparing these to the update formulas. Moreover, as long as the
algorithm does not reach a stationary point there will always be a
strict increase in the objective function $F$. In the third place,
all the $W^n$ and $H^n$ evolve in a compact set. For the $H^n$
this is trivial, since they are nonnegative row stochastic
matrices. For the $W^n$ this follows from~(\ref{eq:wbar2}), since
$W^{n+1}_{il}\leq \sum_j V_{ij}$ (starting the algorithm with
matrices that have strictly positive elements ensures that all
$W^n$ and $H^n$ have strictly positive elements). A detailed
account of the properties of the algorithm is deferred to another
publication.

\end{document}